\documentclass[12pt,a4paper]{article}

\newtheorem{theorem}{Theorem}[section]
\newtheorem{lemma}{Lemma}[section]
\newtheorem{proposition}{Proposition}[section]

\usepackage[dvips]{graphicx}
\usepackage{epsfig}
\usepackage{epsf}

\usepackage{amssymb,amsmath,cite}

\newcommand\sgn{\mathop{\rm sgn}\nolimits}

\begin{document}
\begin{center}\bfseries\Large
An LMI Approach to Stability Analysis of PWM DC-DC Buck Converters
\end{center}

\centerline{\scshape  Alexander N. Churilov}
 \medskip
{\small
\centerline{Department of Computer Science}
\centerline{St. Petersburg Marine Technical University}
\centerline{Lotsmanskaya Str. 3, 190008,
St. Petersburg, Russia}
\centerline{a\_churilov@mail.ru} }
\medskip

\begin{abstract}

{\em The paper considers a DC-DC buck power converter employing
pulse width modulation and voltage feedback control. Global
asymptotical stability of a periodical operating mode is examined.
It is shown that the stability analysis can be reduced to study of
a feasibility problem of semidefinite programming.
}

\end{abstract}

\section{Introduction}

DC-DC buck (step down) power converters serve to convert a direct
current (DC) voltage level to a lower DC voltage level
\cite{Eri99}. One of the most popular types of power converters
uses pulse width modulation (PWM). A pulse width modulator
produces a train of square pulses with a constant frequency and
with a variable duration. By controlling the duration one can
control the output voltage. Last decade DC-DC PWM converters
attracted much attention not only from engineers, but also from
mathematicians and physicists. Their study lies in the mainstream
of common interest to discontinuous and hybrid dynamical systems.

PWM converters may be described by nonlinear
functional-differential equations, their operating modes
correspond to periodic solutions. Many publications were devoted
to the existence of periodic or quasi-periodic solutions, to the
analysis of local stability of such solutions, to bifurcations and
ways of chaotization (see, e.g.,
\cite{FO96,BBC98,YBOY98,ZM03,Tse03}). However,  most of
investigations do not concern global behavior of solutions, i.e.,
a behavior for all possible initial conditions.

There are two basic classes of mathematical models used
for analysis and design of
PWM power converters: discrete time models
\cite{LIYT79,VEK86,HB91,FA02}
and models based on averaging in time of the modulator's
output or of the state space vector
\cite{Sir89b,KBBL90,SNLV91,BL96,TMN04}.
The applicability of the averaging technique suggests
that the switching frequency is high
when compared to the frequencies of
continuous-time signals.

In this paper we use another averaging approach, initiated by the
work \cite{Gel82} and further developed in \cite{GC98}. It is
important to emphasize that our approach is not asymptotical or
approximate, unlike the other averaging technique. However, it
gives conservative estimates if the switching frequency is not
high enough. Our method is based on ideas of the absolute
stability theory \cite{YLG04} and leads to a system of linear
matrix inequalities (LMI) \cite{BGFB94}. This work is a natural
extension of the LMI technique proposed in \cite{CG03}. Some other
methods for stability study of switched power converters via LMIs
can be found in \cite{Rub03,AJKM03}.

The main purpose of this paper is to provide a new computationally
tractable procedure for stability investigations of PWM buck
converters. It can be easily implemented with the help of a
recently developed software for MATLAB and Scilab modelling
systems, such as LMILab \cite{BCPS05}, LMITOOL \cite{END95} or
SeDuMi solver \cite{Stu99} interfaced with YALMIP \cite{Lof04} or
\emph{cvx} \cite{GBY06}.

\section{Preliminaries}

\subsection{Power Stage}\label{S2}

The part of a converter without a control loop is called a power
stage. A conventional power stage of a buck converter is shown
in Fig.~1 (see, e.g., \cite{Eri99}).
Here  $V_s$ is an input voltage, $S$ is an on-off switch,
$L$ is an inductance,  $C_0$ is a capacitance
and $R$ is a load resistance.
Let $U$ be a capacitor voltage (which is the output voltage at the
same time) and $i_L$ be an inductor current.
Then we get
\begin{equation}                         \label{5}
\frac{dU}{dt}=\frac 1{C_0}\, i_L -\frac 1{RC_0}\, U, \qquad
\frac{di_L}{dt}=-\frac 1L\, U+\frac 1L\, V_s,
\end{equation}
when the switch is closed and
\begin{equation}                        \label{6}
\frac{dU}{dt}=\frac 1{C_0}\, i_L -\frac 1{RC_0}\, U, \qquad
\frac{di_L}{dt}=-\frac 1L\, U,
\end{equation}
when the switch is open.

The switch is on with some constant switching period $T$.
Let $\tau$ be an on-time of the switch, i.e., a time
when the switch is closed during the period.
The value $D=\tau/T$ is called duty ratio.

The above described power stage consists of three parts --- a
voltage source, a low-pass filter and a load. The low-pass filter
is designed to remove higher harmonics and pass only a constant
component of the output voltage. However, in fact the output
voltage has a small ripple. Neglecting this ripple, one obtains
$$
U \approx \frac{\tau}T\,V_s=DV_s, \qquad 0\le D\le 1,
$$
so the input voltage steps down.

\begin{figure}
\begin{center}
\includegraphics{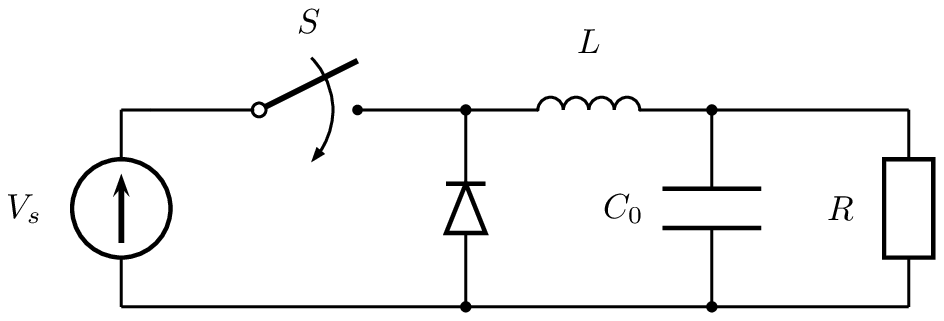}
\end{center}
\caption{Equivalent circuit for a basic power stage}
\end{figure}

In general, the power stage of a buck converter can be described
with a system of equations
\begin{equation}              \label{10}
\frac{dx_p}{dt}= A_px_p+ B_pf,\qquad U= C_px_p
\end{equation}
with
\begin{equation}           \label{11}
\begin{aligned}
f(t)&=\begin{cases}
1, \quad & nT\le t<nT+\tau,
\\
0, \quad & nT+\tau \le t<(n+1)T,
\end{cases}
\qquad
 n&=0,1,\ldots
\end{aligned}
\end{equation}
The signal
$f(t)$ is a train of square single-sign pulses of width
$\tau$ and of period $T$.
The matrix $A_p$ is square, $B_p$ is a column and
$C_p$ is a row.
The matrix $A_p$ is Hurwitz stable, i.e., all its eigenvalues
have negative real parts.

E.g., for equations (\ref{5}), (\ref{6}) we have
\begin{equation*}                         
A_p=\begin{bmatrix} 0 & -1/L \\ 1/C_0 & -1/RC_0
\end{bmatrix},
\qquad
B_p=\begin{bmatrix} V_s/L \\ 0
\end{bmatrix},
\qquad C_p=[0\; 1], \qquad x_p= \begin{bmatrix} i_L \\ U
\end{bmatrix}.
\end{equation*}

Equations (\ref{10}), (\ref{11}) are the equations of an open
loop system (without a control feedback).
In practice, the pulse width $\tau$ is varied to control the
output voltage $U$.
Notice that equations (\ref{10}), (\ref{11}) apply to so called
continuous conduction mode, which suggests that the inductance
$L$ is sufficiently large.

\subsection{Pulse-Width Modulation}\label{S3}

In DC--DC converters a single sign modulation is commonly
employed, i.e. a modulator's output takes only two values, 0 and
1.

In this paper we will limit ourselves to a specific type of
pulse-width modulation called natural or running modulation. The
output of a modulator is
\begin{equation}           \label{13}
f(t)=\begin{cases} 1, \quad & nT\le t<nT+\tau_n,
\\
0, \quad & nT+\tau_n \le t<(n+1)T,
\end{cases}
\qquad
 n=0,1,\ldots
\end{equation}
An input signal $\sigma(t)$
is compared with a sawtooth oscillating signal (a ramp)
$$
\sigma_r(t)=\sigma_1+\sigma_*(t-nT)/T,
\quad  nT\le t< (n+1)T.
$$
Here $\sigma_*$, $\sigma_1$ are given positive parameters.
A switching instant is determined from the relationship
$\sigma(t)=\sigma_r(t)$ (see Fig.~2).

\begin{figure}
\begin{center}
\includegraphics{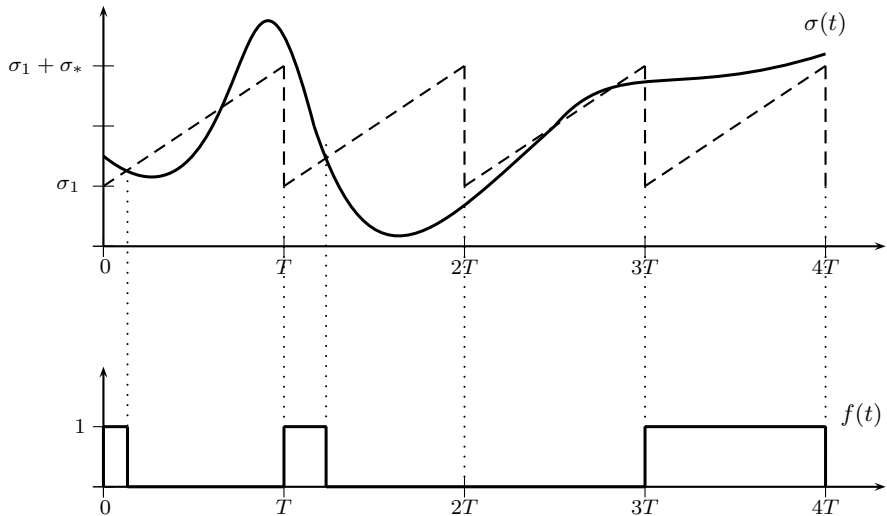}
\end{center}
\caption{A pulse-width modulation}
\end{figure}

More accurately, the modulation law can be described as
follows (see, e.g., \cite{ZM03}).
If $\sigma(nT)\le \sigma_1$, then $\tau_n=0$.
If
\begin{equation}                         \label{15}
\sigma(nT+\tau)>\sigma_1+\sigma_*\tau/T, \quad 0\le\tau\le T,
\end{equation}
then $\tau_n=T$. In all other cases $\tau_n$ is the smallest value
$\tau\in[0,\,T]$ to satisfy the equation
\begin{equation}            \label{16}
\sigma(nT+\tau)=\sigma_1+\sigma_*\tau/T.
\end{equation}
All the other roots of equation (\ref{16}) are ignored.

Such type of modulation is sometimes called modulation with a
latch. In this case every switching interval contains at most one
pulse, so such unpleasant effect as chattering cannot take place.

Besides the above scheme, where the
trailing edge of a pulse is modulated,
there are other modulation laws with modulation of the
front edge or of both edges \cite{LM98}.

\subsection{Control Loop}\label{S4}

Here we consider a buck converter
with a voltage feedback control.
In most cases a closed loop system with a voltage mode control
can be described by the equations
\begin{equation}                   \label{17}
\frac{dx}{dt}=Ax+Bf+q,\quad \sigma=Cx+\psi, \quad f=M\sigma.
\end{equation}
Here $A$ is a constant square matrix, $B$ and $q$ are columns, $C$
is a row, and $\psi$ is a scalar. An operator $M$ describes a
pulse modulator, so functions $\sigma(t)$ and $f(t)$ are the input
and the output of the modulator, respectively. Suppose $A$ to be
Hurwitz stable, i.e., all of its eigenvalues lie in the open left
half-plane.

In the simplest case the control signal is
defined by the formula
\begin{equation}                   \label{18}
\sigma(t)=a(V_{ref}-U(t)).
\end{equation}
Here $V_{ref}$ is a constant reference signal and $a$ is a gain.
Combine (\ref{10}) and (\ref{18}) to obtain (\ref{17}) with
\begin{equation*}                       
A=A_p,\qquad B=B_p, \qquad C=-aC_p,\qquad \psi=aV_{ref}, \qquad
q=0.
\end{equation*}

In a more general case  two additional linear circuits are added
to the control scheme, as shown in Fig.~3 (see, e.g., \cite{Eri99}).
\begin{figure}
\begin{center}
\includegraphics{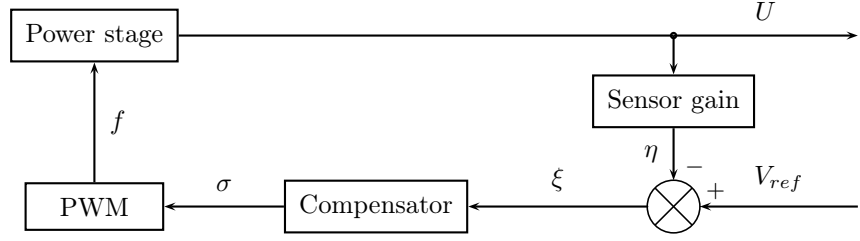}
\end{center}
\caption{A closed loop system}
\end{figure}

Let state space realizations of the compensator and of the sensor gain be
\begin{equation}                       \label{20}
\frac{dx_c}{dt}=A_cx_c+B_c\xi, \qquad \sigma=C_cx_c+D_c \xi
\end{equation}
and
\begin{equation}                       \label{21}
\frac{dx_s}{dt}=A_sx_s+B_s U, \qquad \eta=C_sx_s+D_s U,
\end{equation}
respectively, and
\begin{equation}                       \label{22}
\xi=V_{ref}-\eta.
\end{equation}

Combining (\ref{20}), (\ref{21}) and (\ref{22})
one obtains (\ref{17}) with
\begin{equation*}                       
\begin{aligned}
A&=\begin{bmatrix} A_p & 0  & 0\\
B_sC_p & A_s & 0 \\
-B_cD_sC_p & -B_cC_s & A_c
\end{bmatrix},
\qquad B=\begin{bmatrix} B_p \\ 0 \\ 0\end{bmatrix}, \qquad
q=\begin{bmatrix}0\\ 0\\ B_cV_{ref}  \end{bmatrix},
\\[4mm]
C&=[-D_c D_s C_p, \quad -D_cC_s,\quad C_c], \qquad
\psi=D_cV_{ref}.
\end{aligned}
\end{equation*}

With the help of the change of variables $\tilde x=x+A^{-1}q$,
equations (\ref{17}) are transformed into the equations
\begin{equation}                   \label{17a}
\frac{d\tilde x}{dt}=A\tilde x+Bf,\qquad \sigma=C\tilde x+\psi,
\qquad f=M\sigma
\end{equation}
with $\psi=\psi_0-CA^{-1}q$.

\section{Operating Modes}\label{S5}

Operating modes of a converter correspond to periodical solutions
of the closed loop system. As a rule, such solution is
either $T$-periodical, or has a period multiple to $T$.

The methods of finding periodical modes of a PWM system are well
known. By straightforward calculations, one can obtain an explicit
formula for the response of the linear part of a system to a
square pulse signal. Taking into account periodicity one comes to
a system of transcendental equations with respect to a pulse
duration, which can be solved by numerical methods
(see, e.g., \cite{GC98,ZM03}).

Another way to find a periodical solution is to apply the
harmonic balance method \cite{SNLV91}.
Certainly, in this case the linear
part of the system needs to have good filtering properties.

In this work we will limit our considerations to
$T$-periodical solutions, i.e., to solutions whose period
coincides with the switching period.
In addition, we require these solutions to be \textit{unsaturated},
i.e., the cases $f(t)\equiv 0$ or $f(t)\equiv 1$ are excluded.

Let us formulate necessary and sufficient conditions for the
existence of a $T$-periodical solution of (\ref{17a}). Define for
any pair $(\tau,\, t)$, $0\le\tau\le T$, $0\le t\le T$, a function
\begin{equation}                              \label{17d}
\hat \sigma(\tau, t)=Ce^{At}\hat x(\tau)+\psi+
C\left(e^{At}-I\right) A^{-1}B
\end{equation}
with
\begin{equation}                              \label{17e}
\hat x(\tau)=-\left( I-e^{-AT}\right)^{-1} \left(I-
e^{-A\tau}\right) A^{-1}B.
\end{equation}
Here $I$ is the identity matrix.

\begin{proposition}                                     \label{T0}
System (\ref{17a}) has an unsaturated $T$-periodical solution with
a pulse duration $\tau$ if and only if the following conditions
are satisfied:
\begin{equation*}
0<\tau<T;\qquad \hat\sigma(\tau,\tau)=\Phi(\tau); \qquad
\hat\sigma(\tau,t)>\Phi(t)\quad \mbox{for all}\quad 0\le t<\tau.
\end{equation*}
Here $\Phi(t)=\sigma_1 + \sigma_* t/T$.

For this solution
\begin{align*}
\sigma(t) &= \hat\sigma(\tau,t), \qquad 0\le t\le\tau,
\\
\sigma(t) &=
\hat\sigma(\tau,t)+C\left[I-e^{A(t-\tau)}\right]A^{-1}B, \qquad
\tau\le t\le T.
\end{align*}
\end{proposition}
\bigskip

Proposition~\ref{T0} is proved readily by a direct computation.
The next theorem provides computationally tractable sufficient
conditions for the existence of an unsaturated $T$-periodical
solution.

\begin{theorem}                            \label{T1}
Consider system (\ref{17a}). Suppose that the matrix $A$ is
Hurwitz stable and the inequality
\begin{equation}                             \label{25}
\sigma_1<\psi<\sigma_1+\sigma_*+CA^{-1}B
\end{equation}
is valid. Let there exist a number $\varepsilon>0$ and a
symmetrical matrix $P$ such that the matrix inequalities
\begin{equation}                             \label{26}
\begin{aligned}
P(A^\top+\varepsilon I)+(A+\varepsilon I)P&\preceq -\frac
1{2\varepsilon}\, BB^\top,
\\
P\succ 0, \qquad CAPA^\top C^\top&< \gamma^2
\end{aligned}
\end{equation}
are satisfied.%
\footnote{The character $\prec$ ($\preceq$) denotes negative
definiteness (negative semi-definiteness).
The character $\succ$ ($\succeq$) denotes positive
definiteness (positive semi-definiteness). The symbol $\top$ denotes
transpose.}
Here $\gamma=\sigma_*/T -\min\{0,\,CB\}$.
Then system (\ref{17a}) has an unsaturated $T$-periodical mode.
\end{theorem}

The first inequality (\ref{26}) is not linear in the assembly of
variables $P$ and $\varepsilon$ (it is bilinear in these
variables). However, if $\varepsilon$ is fixed, we get an LMI in
$P$. So, (\ref{26}) can be checked by varying $\varepsilon$ in a
loop. Obviously, the first inequality of (\ref{26}) implies
that the eigenvalues of the matrix $A+\varepsilon I$ have
non-positive real parts. Thus
$\varepsilon$ can be varied from zero to
$-\mathrm{Re}\,\lambda_1$, where $\lambda_1$ is an eigenvalue of
$A$ with the minimal absolute value of a real part.

\section{Stability Conditions}\label{S6}

The following statement is the main result of this paper.

\begin{theorem}                            \label{T2}
Assume that there exists a $T$-periodical solution
$x^0(t)$, $\sigma^0(t)$ of (\ref{17}) and an
estimate
\begin{equation}                            \label{27}
\left| \frac{d\sigma^0(t)}{dt} \right| \le L_1
\end{equation}
is valid, where $L_1$ is a positive constant.
Let $\tau^0$ be a pulse duration for this solution.

Let there exist a symmetrical $m\times m$ matrix $H$ and
scalars $\varepsilon$, $\nu$
such that a system of linear matrix inequalities
\begin{equation}                            \label{29}
L(H)\prec R(\varepsilon,\nu),
\quad H\succ 0, \quad \nu>0, \quad \varepsilon>0
\end{equation}
is satisfied.
Here $L$, $R$ are square symmetrical matrices,
$$
L(H)=\left[
\begin{array}{lccc}
HA+A^\top H & HAB & HB & 0 \\[1mm]
B^\top A^\top H & 0 & 0 & 0 \\[1mm]
B^\top H & 0 & 0 & 0 \\[1mm]
0 & 0 & 0 & 0
\end{array}
\right],
$$
\bigskip
$$
R(\varepsilon,\,\nu)= \left[
\begin{array}{c@{\extracolsep{8mm}}ccc}
0 & 0 &-\frac 12\, C^\top & \kappa_2 A^\top C^\top
\\[1mm]
0 & 3\nu/T^2 & 0 & -\kappa_2\kappa_1
\\[1mm]
-\frac 12\, C & 0 & \sigma_*-TL_1-\nu-\varepsilon-T|\kappa|
& -\kappa_2\kappa
\\[1mm]
\kappa_2 CA & -\kappa_2\kappa_1 & -\kappa_2\kappa & \varepsilon
\end{array}
\right]
$$
with $\kappa=-CB$, $\kappa_1=-CAB$, $\kappa_2=T/\pi$.

Then for any solution $x(t)$ of (\ref{17}) the limit relationships
\begin{equation}                              \label{29a}
\lim_{t\to +\infty}[x(t)-x^0(t)]=0, \qquad
\lim_{n\to\infty} \tau_n=\tau^0
\end{equation}
are satisfied.
Moreover, the solution $x^0(t)$ is stable in the sense of Lyapunov.
\end{theorem}

Inequalities (\ref{29}) can be easily verified with the help of
standard software packages that are suitable for study an LMI
feasibility problem \cite{BGFB94}.

\section{Numerical Example}\label{S7}

Consider a converter with a power stage described by (\ref{5}),
(\ref{6}), where
$R=22\,\Omega$,
$C_0=47\,\mu F$, $L=20\,mH$, $T=400\,\mu s$, $V_s= 20\,V$.
Let a control be described by (\ref{18}) with $a=1$,
$V_{ref}=13.5\,V$, $\sigma_1=4\,V$. The parameter $\sigma_*$
will be chosen later.
Then we have $\psi=13.5$, $q=0$.

Conditions of Theorem~\ref{T0} are satisfied for
\begin{equation}                        \label{30}
\sigma_* \ge 12.83\,V.
\end{equation}
Thus, under condition (\ref{30}) the system has an unsaturated
$T$-periodical mode.

Apply MATLAB to get a response of the power stage to the pulse
signal with a fixed duty ratio $\tau/T$.
The direct modelling gives the following results:
\begin{center}
\begin{tabular}{|c|c|c|c|c|c|}
\hline
$\tau/T$ & 0.1 & 0.3 & 0.5 & 0.7 & 0.9
\\ \hline
$TL_1$ & 0.15 & 0.35 & 0.44 & 0.35 & 0.15
\\ \hline
\end{tabular}
\end{center}
Thus
$TL_1=0.44$ can be chosen for the worst case.

Applying Theorem~\ref{T1} and the standard LMILab package for MATLAB
gives
\begin{equation}                        \label{32}
\sigma_* \ge 17.78\,V.
\end{equation}
Obviously, (\ref{30}) follows from (\ref{32}).
In result, under condition (\ref{32}) the system has a globally
stable $T$-periodical mode.

If we neglect the ripple of the output voltage,
the following approximate formula can be easily obtained:
$$
U^0\approx V_s\,\frac{a V_{ref}-\sigma_1}{\sigma_*+aV_s}.
$$
Choose $\sigma_*=18\,V$, which satisfies (\ref{32}).
Then $U^0=5\,V$, so the converter steps the voltage down
from $20\,V$ to $5\,V$.

\section{Proofs of the Theorems}

\textbf{A Useful Lemma.}
We precede the proof of Theorem~\protect\ref{T1} with a lemma,
which is a version of a statement proved in \cite{GC86}.

\begin{lemma}                          \label{L1}
Suppose that $A$ is a square matrix, $B$ is a column and
$\varepsilon$ is a positive scalar. Then the inequality
\begin{equation}                        \label{33}
H(A+\varepsilon I)+(A^\top +\varepsilon I)H \preceq
-\frac 1{2\varepsilon} HBB^\top H
\end{equation}
is equivalent to the relationship
\begin{equation}                        \label{34}
x^\top H(Ax+Bf) +\varepsilon (x^\top Hx-1) \le 0,
\qquad \forall x, \; \forall f, \; |f|\le 1.
\end{equation}
\end{lemma}

Since
$$
2x^\top HBf-2\varepsilon-\frac 1{2\varepsilon} (x^\top HB)^2f^2 =
-\left(\frac 1{\sqrt{2\varepsilon}}\,x^\top HBf- \sqrt{2\varepsilon}
\right)^2, \qquad \forall x,\;\forall f,
$$
(\ref{34}) follows from  (\ref{33}).

Let  (\ref{34}) be satisfied. Obviously, (\ref{34})  preserves if
one replaces $x$ for a vector $\alpha x$, where $\alpha$ is an
arbitrary scalar. Let us fix arbitrary $x$ and $f$. Then (\ref{34})
implies
$$
\alpha^2 x^\top H(A\varepsilon I) +\alpha x^\top HBf-\varepsilon\le 0,
\qquad \forall \alpha.
$$
The latter relationship is satisfied for all $\alpha$ if and only if
$$
(x^\top HBf)^2 +4\varepsilon x^\top H(A+\varepsilon I) x \le 0.
$$
Since $|f|\le 1$, we come to (\ref{33}).

\textbf{Proof of Theorem~\protect\ref{T1}:}
The proof generally follows the proof of Theorem~6.2 \cite{GC98}.
Let us check that under condition (\ref{25})  system (\ref{17a})
has no saturated periodical solutions (equilibria). Firstly,
suppose that (\ref{17a}) has a solution with $f(t)\equiv 0$, that
is $\tau=0$. From the modulation law this means that $\sigma(0)\le
\sigma_1$. Hence, from (\ref{17a}) one concludes that
$\psi\le\sigma_1$. This contradicts the first inequality
(\ref{25}). Similarly, let $f(t)\equiv 1$, that is $\tau=T$. From
(\ref{17d}), (\ref{17e}) it follows that
$\tilde\sigma(t)=\hat\sigma(T,t)=\psi-CA^{-1}B$. Then (\ref{15}),
(\ref{16}) imply $\psi-CA^{-1}B\ge\sigma_1+\sigma_*$, which
contradicts the second inequality (\ref{25}). Thus system
(\ref{17a}) cannot have an equilibrium.

Let us prove the existence of a $T$-periodical solution under
condition  (\ref{26}).  Define $H=P^{-1}$. Then (\ref{26}) imply
$H\succ 0$, (\ref{33}), and
\begin{align}                               \label{40}
CAH^{-1}A^\top C^\top &<\gamma^2.
\end{align}

Firstly, consider the case, when $CB\le 0$ and $\gamma=\sigma_*/T-CB$.
Let us take an ellipsoid
$$
E=\left\{\, x\;:\; x^\top Hx\le 1\,\right\}
$$
and define a Lyapunov function $V(x)=x^\top Hx$. From Lemma~\ref{L1}
it follows that
$$
\frac{dV}{dt}+2\varepsilon(V-1) \le 0
$$
along the solutions of (\ref{17a}). Thus if a solution starts in
$E$, it remains in $E$ for all subsequent times. Since
$$
\max_{x\in E} CAx=(CAH^{-1}A^\top C^\top)^{1/2},
$$
from inequality (\ref{40}) it follows that
\begin{equation}                               \label{45}
\max_{x\in E} CAx <\sigma_*/T-CB.
\end{equation}

Consider the operator of translation along the trajectories
of (\ref{17a}):
$$
S\;:\; \tilde x(0) \mapsto \tilde x(T).
$$
As it was shown above, $S E\subset E$. Let us prove that $S$ is
continuous on $E$. It is easily seen that $S$ is continuous
provided that the pulse width $\tau_0$ depends continuously on
$\tilde x(0)$. Since $\tau_0$ is the minimal positive root of
equation (\ref{16}) with $n=0$, it suffices to require
$$
\frac{d\sigma(\tau_0-0)}{d\tau} \ne \frac{\sigma_*}T.
$$
The last inequality can be rewritten as
$$
CA\tilde x(\tau_0)+CB \ne \sigma_*/T.
$$
Thus continuity of $S$ follows from (\ref{45}).

We have found that
the operator $S$ is continuous on a closed, bounded and convex set.
Hence, it has a fixed point, which corresponds to a $T$-periodical
solution~\cite{Kra68}.

Let us turn to the case, when $CB>0$ and $\gamma=\sigma_*/T$.
Consider an ellipsoid
$$
E_0=\left\{\, x\;:\; (x+A^{-1}B)^\top H(x+A^{-1}B)\le 1\,\right\}
$$
and a Lyapunov function $V_0(x)=(x+A^{-1}B)^\top H(x+A^{-1}B)$.
Obviously,
$$
\frac{dV_0}{dt}=2(x+A^{-1}B)^\top H[A(x+A^{-1}B)+B(f-1)]
$$
along the solutions of (\ref{17a}).
Arguing as above, one obtains that the ellipsoid $E_0$
is invariant under the translation $S$ and $S$ is continuous
on $E_0$. Hence, it has a fixed point on this ellipsoid.

\textbf{Proof of Theorem~\protect\ref{T2}:}
The proof is based on the proofs of Theorems~3.3, 4.7
\cite{GC98} and follows the averaging scheme proposed by
A.~Kh.~Gelig in \cite{Gel82}.
The idea of the approach is that the state vector $\tilde x(t)$
varies much slower than the modulator's output $f(t)$, so
this output can be averaged in time \cite{And60}.

Let $x(t)$, $\sigma(t)$ be an arbitrary solution of (\ref{17}).
Consider the deviations $x_d(t)=x(t)-x^0(t)$,
$\sigma_d(t)=\sigma(t)-\sigma^0(t)$, $f_d(t)=f(t)-f^0(t)$, where
$f^0=M\sigma^0$.
Then we obtain
\begin{align}                              \label{50}
\frac{dx_d}{dt}&=Ax_d+Bf_d, \qquad \sigma_d=Cx_d,
\\[2mm]
f_d(t) &= \begin{cases}
0, & nT \le t< nT+\tau_n^{\min},
\\
\sgn (\tau_n -\tau^0), \qquad & nT+\tau_n^{\min}\le t< nT+\tau_n^{\max},
\\
0, & nT+\tau_n^{\max} \le t <nT+T.
\end{cases}
\end{align}
Here $\tau_n^{\min}=\min\{\tau_n,\tau^0\}$,
$\tau_n^{\max}=\max\{\tau_n,\tau^0\}$. Recall that $\tau^0$ is
the pulse duration for the given $T$-periodical mode.

Define a sequence of average values
$$
v_n=\frac 1T\int\limits_{nT}^{nT+T}
f_d(t)\,dt=\frac{\tau_n-\tau^0}T
$$
and consider a piecewise constant function
$$
v(t)=v_n, \qquad nT\le t<nT+T.
$$
Unlike the signal $f_d(t)$, which is modulated in width and position,
the pulse signal $v(t)$
is amplitude modulated. However, the average values of both signals
coincide on each sampling interval. One has
$$
f_d(t)=v(t)+\frac{du(t)}{dt},
\qquad
\mbox{where}
\qquad
u(t)=\int\limits_0^t (f_d(s)-v(s))\,ds.
$$
Obviously, $u(nT)=0$ for $n=0,1,\ldots$. The function $u(t)$
may be considered as an averaged error of the replacement of
$f_d(t)$ for $v(t)$. In Chapter~3, \cite{GC98} it was shown that $u(t)$
is much less than $v(t)$, provided that the switching period $T$ is small.
More precisely, $|u(t)|\le T|v(t)|$ for all $t$ and
$$
\int\limits_{nT}^{nT+T} u(t)^2\,dt\le
\frac{T^2}3 \int\limits_{nT}^{nT+T} v(t)^2\,dt, \qquad
n=0,\,1,\ldots.
$$

The following lemma plays the main role in the proof.

\begin{lemma}                       \label{L2}
For any solution of (\ref{17})  and
for any $n$ there exists a number $\tilde t_n$ such that
$nT\le \tilde t_n\le nT+T$ and
\begin{equation}                              \label{60}
0\le \frac{v_n}{\sigma_d(\tilde t_n)} \le  \frac 1{\sigma_*-TL_1}.
\end{equation}
\end{lemma}

Inequalities (\ref{60}) present a version of Lur'e sectoral
constraints, which are common in the absolute stability theory.

\textbf{Proof of Lemma~\ref{L2}:} Since the given $T$-periodic mode is
unsaturated, the relationship
\begin{equation}                              \label{61}
\sigma^0(\tau^0)=\sigma^0(\tau^0+nT)=
\sigma_1+\sigma_*\tau^0/T
\end{equation}
is valid. Consider the following three cases.

Case (i). Suppose that $\tau_n$ is defined as a minimal positive
root of equation (\ref{16}). Subtract
$\sigma^0(nT+\tau_n)=\sigma^0(\tau_n)$ from both sides of equality
(\ref{16}) to obtain
\begin{equation}                              \label{62}
\sigma_d(nT+\tau_n)=\sigma^0(\tau^0)-\sigma^0(\tau_n)+\sigma_* v_n.
\end{equation}
Since $\tau_n=\tau^0+Tv_n$, one can write
\begin{equation}                              \label{63}
\sigma^0(\tau_n)-\sigma^0(\tau^0)=
\sigma^0(\tau^0+Tv_n)-\sigma^0(\tau^0)=Tv_n\frac{d\sigma^0}{dt}(t_s),
\end{equation}
where $t_s$ is some time instant lying between $\tau^0$ and
$\tau_n$. Then (\ref{62}) and (\ref{63}) imply
$$
\frac{v_n}{\sigma_d(nT+\tau_n)}=\frac 1{\sigma_*-T\dfrac{d\sigma^0}{dt}(t_s)},
$$
so
$$
\frac 1{\sigma_*+TL_1}\le \frac{v_n}{\sigma_d(nT+\tau_n)}
\le  \frac 1{\sigma_*-TL_1}.
$$
Hence condition (\ref{60}) is satisfied with $\tilde t_n=nT+\tau_n$.

Case (ii). Suppose that
\begin{equation}                              \label{64}
\sigma(nT)<\sigma_1,
\end{equation}
so $\tau_n=0$ and
$v_n=-\tau^0/T<0$. Let us estimate $v_n$ from below.
In view of (\ref{61}), one has
\begin{equation}                              \label{65}
\sigma_* v_n=-\sigma_*\tau^0/T=\sigma_1-\sigma^0(\tau^0).
\end{equation}
Then  (\ref{64}) and (\ref{65}) imply
\begin{equation}                              \label{66}
\sigma_* v_n > \sigma(nT)-\sigma^0(\tau^0)=\sigma_d(nT)+\sigma^0(nT)-
\sigma^0(\tau^0).
\end{equation}
Write
$$
\sigma^0(0)-\sigma^0(\tau^0)=Tv_n \frac{d\sigma^0}{dt}(t_s),
$$
where $t_s$ is some time,  $0\le t_s\le \tau^0$. Since
$\sigma^0(nT)=\sigma^0(0)$, (\ref{66}) implies
$$
0>v_n\left( \sigma_* - \frac{d\sigma^0}{dt}(t_s) \right)>\sigma_d(nT).
$$
Then we come to (\ref{60}) with $\tilde t_n=nT$.

Case (iii). Suppose that (\ref{15}) is fulfilled, so $\tau_n=T$
and $v_n = 1-\tau^0/T >0$.  Set $\tau=0$ and $\tau=T$ in (\ref{15}).
Then
$$
\sigma(nT)>\sigma_1, \qquad\sigma(nT+T)>\sigma_1+\sigma_*.
$$

Firstly, consider the case when  $\sigma(nT)\ge \sigma_1+\sigma_*$.
Then
$$
\sigma_* v_n= \sigma_1+\sigma_*-\sigma^0(\tau^0)\le
\sigma(nT)-\sigma^0(\tau^0)=\sigma_d(nT)+\sigma^0(0)-\sigma^0(\tau^0).
$$
Arguing as above, we get
$$
0<v_n\left( \sigma_* - \frac{d\sigma^0}{dt}(t_s) \right)<\sigma_d(nT),
$$
so we come to (\ref{60}) with $\tilde t_n=nT$.

At last suppose that  (\ref{15}) is  satisfied and
$\sigma(nT)< \sigma_1+\sigma_*$.  Since
$\sigma(nT+T)> \sigma_1+\sigma_*$, there exists a number
$\tau_s$ such that $0<\tau_s<T$ and $\sigma(nT+\tau_s)=\sigma_1+\sigma_*$.
It is easily seen that (\ref{60}) is valid with $\tilde t_n=nT+\tau_s$.
The proof of Lemma~\ref{L2} is complete.

Now return to the proof of Theorem~\ref{T1}.
Inequality  (\ref{29}) can be rewritten as
\begin{equation}                              \label{67}
L_0(H)\prec R_0(\varepsilon,\,\nu)
\end{equation}
with
$$
L_0(H)=\left[
\begin{array}{lcc}
HA+A^\top H & HAB & HB \\[1mm]
B^\top A^\top H & 0 & 0  \\[1mm]
B^\top H & 0 & 0
\end{array}
\right],
$$
\bigskip
$$
R_0(\varepsilon,\,\nu)= \left[
\begin{array}{ccc}
0 & 0 &-\frac 12\, C^\top
\\[1mm]
0 & 3\nu/T^2 & 0
\\[1mm]
-\frac 12\, C & 0 & \sigma_*-TL_1-\nu-\varepsilon-T|\kappa|
\end{array}
\right]  -
\frac{\kappa_2^2}{\varepsilon}
\begin{bmatrix}
A^{\top}C^{\top} \\ -\kappa_1 \\ -\kappa
\end{bmatrix}
[CA\;  -\kappa_1 \; -\kappa].
$$
Inequality (\ref{67}) has less dimension than (\ref{29}),
however, it is not linear in $\varepsilon$.

Let us make a change of variables $y(t)=x_d(t)-Bu(t)$.
Then (\ref{50}) can be rewritten as
\begin{equation}                              \label{70}
\frac{dy}{dt}=Ay+Bv+ABu, \qquad\sigma=Cy+CBu.
\end{equation}
Consider a Lyapunov function $V(y)=y^\top Hy$, where $H$
satisfies  (\ref{67}). Define a quadratic form
$$
F(y,u,v)=
\begin{bmatrix}
y\\u\\v
\end{bmatrix}^\top
R_0(\varepsilon,\,\nu)
\begin{bmatrix}
y\\u\\v
\end{bmatrix}.
$$
Then (\ref{67}) implies
\begin{equation}                              \label{71}
\frac{dV(y(t))}{dt}\le F(y,u,v) -\delta (\|y\|^2+u^2+v^2)
\end{equation}
along the solutions of (\ref{70}). Here $\delta$ is some
small positive number. From Lemma~2 it follows that
$$
\int\limits_{nT}^{nT+T} F(y(t),u(t),v(t))\,dt \le 0
$$
for all $n$ (see Chapter~3 \cite{GC98} for details). Then (\ref{71})
implies $y(t)\to 0$, $v(t)\to 0$ as $t\to +\infty$.
Thus we come to (\ref{29a}). The Lyapunov stability can be proved
as in  \cite{CG03}.

\section {Conclusion}

We propose an LMI based approach to
finding conditions for global stability of an operating mode
of a PWM power converter.
The LMI technique becomes a powerful tool for
this type of problems
when combined with well-known methods of
calculation of periodical solutions for PWM systems.
The LMI approach leads to computationally tractable criteria,
which can be easily implemented with the help of standard
modelling software.

\section*{Acknowledgment}

The work was supported in part by the
Russian Foundation for Basic Research, project
05-01-00290-a.

\end{document}